\documentclass[11pt]{amsart}

\usepackage{amssymb}
\usepackage{amsmath}
\usepackage{amssymb}
\usepackage{graphicx}
\usepackage{amsmath,mathrsfs}
\usepackage[all]{xy}

\theoremstyle{plain}

\newtheorem*{theorem*}{Theorem*}

\theoremstyle{remark}

\newcommand{\Set}{\mathsf{Set}}

\newcommand{\T}{\mathcal{T}}
\newcommand{\F}{\mathcal{F}}
\newcommand{\A}{\mathcal{A}}
\newcommand{\C}{\mathcal{C}}
\newcommand{\E}{\mathcal{E}}
\newcommand{\M}{\mathcal{M}}
\newcommand{\B}{\mathcal{B}}
\newcommand{\W}{\mathcal{W}}

\newcommand{\ob}{\mathrm{ob}}

\newcommand{\iso}{\mathrm{Iso}}
\newcommand{\cart}{\mathrm{Cart}}

\newdir{|>}{%
!/4.5pt/@{|}*:(1,-.2)@^{>}*:(1,+.2)@_{>}}

\long\def\symbolfootnote[#1]#2{\begingroup%
\def\thefootnote{\fnsymbol{footnote}}\footnote[#1]{#2}\endgroup}

\begin{document}

\title{Factorization, Fibration and Torsion}

\author{Ji\v r\'{i} Rosick\'{y}$^*$}

\address{Department of Algebra and Geometry\\
Masaryk University\\
Jan\'{a}\v{c}kove n\'{a}m. 2a, 66295 Brno, Czech Republic}
\email{rosicky@math.muni.cz}

\author{Walter Tholen$^\dag$}

\address{Department of Mathematics and Statistics\\
  York University\\
 Toronto, ON M3J 1P3, Canada
}
\email{tholen@mathstat.yorku.ca}

\date{\today}

\begin{abstract}
A simple definition of torsion theory is presented, as a
factorization system with both classes satisfying the 3--for--2
property. Comparisons with the traditional notion are given, as well
as connections with the notions of fibration and of weak
factorization system, as used in abstract homotopy theory.

\end{abstract}


\maketitle
\begin{center}
    \textsc{\small{Dedicated to the memory of Saunders Mac Lane} }
\end{center}

\section{Introduction}

\symbolfootnote[0]{$^\dag$ Partial financial assistance by NSERC
is gratefully acknowledged.}

\symbolfootnote[0]{$^{*}$ Supported by the Grant Agency of the Czech
Republic under grant no. 201/06/0664. The hospitality of York
University is gratefully acknowledged.}

That full reflective subcategories may be characterized by certain
factorization systems is well known, thanks to the works of Ringel
[Ri] and Cassidy, H\'{e}bert and Kelly [CHK]. While the former paper
treats the characterization in the context of the Galois
correspondence that leads to the definition of weak factorization
systems (as given in [AHRT]), the latter paper carefully analyzes
construction methods for the factorizations in question. To be more
specific, following [CHK], we call a factorization system ($\E, \M)$
\emph{reflective} if $\E$ satisfies the cancellation property that
$g$ and $gf$ in $\E$ force $f$ to be in $\E$; actually, $\E$ must
then have what homotopy theorists call the 3-for-2 property. When
there is a certain one--step procedure for constructing such
factorizations from a given reflective subcategory, the system is
called \emph{simple}. Following a pointer given to the second author
by Andr\'{e} Joyal, in this paper we characterize simple reflective
factorization systems of a category $\C$ in terms of generalized
fibrations $P: \C \to \B$: they are all of the form $\E = \{
\text{morphisms inverted by } P\}$, $\M = \{ P-\text{cartesian
morphisms}\}$ (see Theorem 3.9). In preparation for the theorem, we
not only carefully review some needed facts on factorization
systems, but characterize them also within the realm of weak
factorization systems (Prop. 2.3), using a somewhat hidden result of
[Ri], and we frequently allude to the use of weak factorization
systems in the context of Quillen model categories. Furthermore, we
have included a new result for many types of categories, including
extensive categories as well as additive categories, namely that
$(\{\text{coproduct injections}\}, \{ \text{split epimorphisms}\})$
form always a weak factorization system (Theorem 2.7), which is
somewhat surprising since both classes appear to be small.

The main point of the paper, however, is to present an easy
definition of torsion theory that simplifies the definition given by
Cassidy, H\'ebert and Kelly [CHK]. Hence, here a torsion theory in
any category is simply a factorization system $(\E, \M)$ that is
both reflective and coreflective, so that both $\E$ and $\M$ have
the 3--for--2 property. At least in pointed categories with kernels
and cokernels, such that every morphism factors into a cokernel
followed by a morphism with trivial kernel, and dually, our torsion
theories determine a pair of subcategories with the properties
typically expected from a pair of subcategories of ``torsion''
objects and of ``torsion--free'' objects, at least when the system
$(\E, \M)$ is simple (Theorem 4.10). We present a precise
characterization of ``standard'' torsion theories (given by pairs of
full subcategories) in terms of our more general notion in Theorem
5.2, under the hypothesis that the ambient category is homological
(in the sense of [BB]), such that every morphism factors into a
kernel preceded by a morphism with trivial cokernel. At least all
additive categories which are both regular and coregular (in the
sense of Barr [Ba]) have that property.

We have dedicated this paper to the memory of Saunders Mac Lane,
whose pioneering papers entitled ``Groups, categories and duality''
(Bulletin of the National Academy of Sciences USA 34(1948) 263-267)
and ``Duality for groups'' (Bulletin of the American Mathematical
Society 56 (1950) 485-516) were the first to not only introduce
fundamental constructions like direct products and coproducts in
terms of their universal mapping properties, but to also present a
forerunner to the modern notion of factorization system, an
equivalent version of which made its first appearance in John
Isbell's paper ``Some remarks concerning categories and subspaces''
(Canadian Journal of Mathematics 9 (1957) 563-577), but which became
widely popularized only through Peter Freyd's and Max Kelly's paper
on ``Categories of continuous functors, I'' (Journal of Pure and
Applied Algebra 2 (1972) 169-191).

Some of the results contained in this paper were presented by the
second author at a special commemorative session on the works of
Samuel Eilenberg and Saunders Mac Lane during the International
Conference on Category Theory, held at White Point (Nova Scotia,
Canada) in June 2006.

\emph{Acknowledgement}: The authors thank George Janelidze for many
helpful comments on an earlier version of the paper, especially for
communicating to them the current proof of Theorem 2.7 which
substantially improves and simplifies their earlier argumentation.

\section{Weak factorization systems and factorization systems}

\subsection* {2.1} For morphisms $e$ and $m$ in a category $\C$ one
writes
\begin{equation*}
e \; \square \; m \;\;\; (e \bot m)
\end{equation*}
if, for every commutative solid-arrow diagram
\begin{equation*}
  \xymatrix{
    . \ar[r] \ar[d]_e & . \ar[d]^m \\
    . \ar[r] \ar@{-->}[ur]^d & .
  }
\end{equation*}
one finds a (unique) arrow $d$ making both emerging triangles
commutative. For classes $\E$ and $\M$ of morphisms in $\C$ one
writes
\begin{equation*}
\begin{array}{cc}
  \E^{\square} = \{m \:|\; \forall e \in\E : e\;\square\;m\}, & ^{\square}\M = \{e\;|\; \forall m \in \M : e \;\square\; m\}, \\
  \E^{\bot} = \{m \;|\; \forall e\in \E :e \bot m\}, & ^{\bot}\M = \{e \;|\; \forall m \in \M : e \bot m\}. \\
\end{array}
\end{equation*}
Recall that $(\E, \M)$ is a \textit{weak factorization system
(wfs)} if
\begin{enumerate}
    \item $\C = \M\cdot\E$
    \item $\E =$ $^{\square}\M$ \textit{and} $\M = \E^{\square}$,
\end{enumerate}
and it is a \textit{factorization system (fs)} if (1) holds and
\begin{enumerate}
    \item[(2*)] $\E =$ $ ^{\bot}\M$ \textit{and} $\M = \E^{\bot}$
\end{enumerate}
It is well known that, in the presence of (1) , condition (2) may
be replaced by
\begin{enumerate}
    \item[(2a)] $\E\square\M$ (that is: $e\;\square\; m$ for all
    $e \in \E$ and $m \in \M$), and
    \item[(2b)] \textit{$\E$ and $\M$ are closed under retracts in
    $\C^2 (=\C^{\{\xymatrix@1@=8pt{\cdot \ar[r] & \cdot}\}})$},
\end{enumerate}
and (2b) may be formally weakened even further to
\begin{enumerate}
    \item[(2b1)] \textit{if $gf \in \E$ with $g$ split mono, then
    $f \in \E$,} and
    \item[(2b2)] \textit{if $gf\in \M$ with $f$ split epi, then $g \in \M$} (see [AHRT]).
\end{enumerate}
Likewise, in the presence of (1), condition (2*) may be replaced
by
\begin{enumerate}
    \item [(2*a)] $\E \bot \M$, and
    \item [(2*b)] \textit{ $\E$ and $\M$ are closed under
    isomorphisms in $\C^2$.}
\end{enumerate}
\subsection*{2.2} \textit{Every factorization system is a wfs} (see
[AHS], [AHRT]), and for every wfs $(\E, \M)$ one has $\E \cap \M
=$ Iso$\C$, $\E$ and $\M$ are closed under composition, $\E$ is
stable under pushout and closed under coproducts, and $\M$ has the
dual properties. For a factorization system $(\E, \M)$, the class
$\E$ is actually closed under every type of colimit and satisfies
the cancellation property
\begin{enumerate}
    \item[(3)] \textit{ if $gf \in \E$ and $f \in \E$, then $g \in
    \E$.}
\end{enumerate}
Using an observation by Ringel [Ri] (see also [T, Lemma 7.1]) we
show that each of these additional properties characterizes a wfs as
an fs.
\subsection*{2.3 Proposition.} \emph{Let $(\E, \M)$ be a wfs of a
category $\C$ with cokernelpairs of morphisms in $\E$. Then the
following conditions are equivalent:}
\begin{enumerate}
    \item[(i)] \emph{$(\E, \M)$ is a factorization system;}
    \item[(ii)] \emph{$\E$ is closed under any type of colimit (in the morphism category of $\C$);}
    \item[(iii)] \emph{for every $e:A \to B$ in $\E$ the canonical
    morphism $e':B+_AB \to B$ lies also in $\E$ (where $B+_AB$ is
    the codomain of the cokernelpair of $e$);}
    \item[(iv)]\emph{$\E$ satisfies condition} (3)\emph{;}
    \item[(v)] \emph{if $gf =1$ with $f \in \E$, then $g \in \E$.}
\end{enumerate}

\begin{proof}
 $(i)\implies(ii)$ and $(i)\implies(iv)$ are well known (see
2.2), and $(iv)\implies(v)$ is trivial. For $(ii) \implies(iii)$
consider the diagram
\begin{equation*}
  \xymatrix{
    A \ar[r]^e \ar[d]_e & B \ar@<2pt>[r]^(0.35){p_1} \ar @<-2pt>[r]_(0.35){p_2} \ar[d]^1 & B+_AB \ar[d]^{e'}\\
    B \ar[r]^1 & B \ar@<2pt>[r]^{1} \ar @<-2pt>[r]_{1} & B
  }
\end{equation*}
where both rows represent cokernelpairs. Since the connecting
vertical arrows $e$ and $1$ lie in $\E$, $e'$ lies also in $\E$,
by hypothesis. For $(v) \implies (iii)$ observe that, since $\E$
is stable under pushout, one has $e'p_1 = 1$ with $p_1 \in \E$, so
that $e' \in \E$ follows. Finally, for $(iii)\implies (i)$,
consider the diagram
\begin{equation*}
  \xymatrix{
    A \ar[r]^u \ar[d]_e & C \ar[d]^m\\
    B \ar[r]_v \ar@<2pt>[ur]^{s} \ar @<-2pt>[ur]_{t}& D
  }
\end{equation*}
with $e \in \E$, $m \in \M$, $se = te = u$ and $ms = mt = v$. The
morphism $r:B+_AB \to C$ with $rp_1 = s$ and $rp_2 = t$ makes
\begin{equation*}
  \xymatrix{
    B+_AB \ar[r]^(0.7)r \ar[d]_{e'} & C \ar[d]^m\\
    B \ar[r]_v & D
  }
\end{equation*}
commute. Hence, by hypothesis, one obtains $w:B \to C$ with $we' =
r$, and
\begin{equation*}
    s = rp_1 = we'p_1 = w = we'p_2= rp_2 = t
\end{equation*}
follows, as desired.
\end{proof}

 Dualizing (part of) the Theorem we
obtain:
\subsection*{2.4 Corollary} \emph{In a category with kernelpairs,
$(\E, \M)$ is an fs if, and only if, it is a wfs and satisfies the
condition}:

$(v^{\text{op}})$ \emph{if $gf = 1$ with $g \in \M$, then $f \in
\M$.}

\subsection*{2.5} If (Epi, Mono) in $\Set$ is the prototype of fs,
then (Mono, Epi) in $\Set$ is the prototype of wfs. But the latter
claim actually disguises a simple general fact which does not seem
to have been stated clearly in the literature yet. In conjunction
with two very special features of $\Set$, namely that 1. every
monomorphism is a coproduct injection and 2. every epimorphism
splits (=Axiom of Choice), the following Proposition and Theorem
give, \emph{inter alia},  the (Mono, Epi) system:

\subsection*{2.6 Proposition.} \emph{In a category with binary
coproducts}, ($^{\square}$SplitEpi, SplitEpi) \emph{is a wfs, and
a morphism $f: A \to B$ lies in} $^{\square}$SplitEpi \emph{if,
and only if, there is some $k:B \to A+B$ with $kf = i:A \to A+B$
the first coproduct injection, and with $<f,1_B>k = 1_B$; in
particular, every coproduct injection lies in }
$^{\square}$SplitEpi.

\begin{proof}
 Every morphism $f: A \to B$ factors as $pi = f$, and the
co-graph $p:= <f,1_B>:A+B \to B$ is a split epimorphism; moreover,
split epimorphisms satisfy condition (2b2) trivially. It now
suffices to prove the given characterization of morphisms in
$^{\square}$SplitEpi, since it shows in particular that coproduct
injections are in $^{\square}$SplitEpi (simply take $k$ to be a
coproduct injection), and since $^{\square}$SplitEpi (like any
class $^{\square}\M$) satisfies (2b1). Given $f \in
^{\square}$SplitEpi one obtains $k$ from $f \square p$:
\begin{equation*}
  \xymatrix{
    A \ar[r]^(0.4)i \ar[d]_{f} & A+B \ar[d]^p\\
    B \ar@{=}[r] \ar@{-->}[ur]^k & B
  }
\end{equation*}
Conversely, having $k$ with $kf = i$ and $pk = 1_B$, consider the
diagram
\begin{equation*}
  \xymatrix{
    A \ar[r]^u \ar[d]_{f} & X \ar @<2pt>[d]^r\\
    B \ar[r]^v & Y \ar@<2pt>[u]^t
  }
\end{equation*}
with $ru = vf$ and $rt = 1_Y$. Then $s:= <u, tv>:A+B \to X$
satisfies
\begin{equation*}
    rsi = ru = vf = vpi,\;\;\; rsj = rtv = v = vpj,
\end{equation*}
with $j$ the second coproduct injection, so that $rs = vp$. Hence,
$d:= sk: B \to X$ satisfies
\begin{equation*}
    df = skf = si = u, \;\;\; rd = rsk = vpk = v,
\end{equation*}
as desired.
\end{proof}

In many important types of categories, the class $^{\Box}$SplitEpi
is remarkably small:

\subsection*{2.7 Theorem} \emph{Let $\C$ be a category with binary coproducts,
and $\mathrm{Sum}$ be the class of all coproduct injections. If
$\mathrm{Sum}$ is stable under pullback in $\C$, or if $\C$ is
pointed and $\mathrm{Sum}$ contains all split monomorphisms, then
$(\mathrm{Sum}, \mathrm{SplitEpi})$ is a wfs in $\C$. The hypotheses
on $\C$ are particularly satisfied when $\C$ is extensive (in the
sense of $\mathrm{[CLW]}$) or just Boolean (in the sense of
$\mathrm{[M]}$), or when $\C$ is an additive category with finite
coproducts.}

\begin{proof}
It suffices to prove that $f:A \to B$ in $^{\square}$SplitEpi is a
coproduct injection. With the (split) monomorphism $k$ as in 1.6,
consider the diagram
\begin{equation*}
  \xymatrix{
    A \ar[r]^f \ar[d]_{1_A} & B \ar[d]^{1_B} \ar[r]^{1_B} & B \ar[d]^k\\
    A \ar[r]_f & B \ar[r]_k & A+B
  }
\end{equation*}
which is composed of two trivial pullback diagrams. By hypothesis,
since $kf$ is a coproduct injection, its pullback $f$ is also one.

If $\C$ is pointed, the morphism $f:A \to B$ in $^{\square}$SplitEpi
is a split monomorphism (since $<1_A,0>kf = <1_A,0>i = 1_A$), and as
such it is a coproduct injection, by hypothesis.
\end{proof}

For the sake of completeness we mention another well-known general
reason for (Mono, Epi) being a wfs in $\Set$:

\subsection*{2.8 Proposition} [AHRT] \emph{In every category with
binary products and enough injectives,} (Mono, Mono$^{\square})$
\emph{is a wfs.}
\begin{flushright}
$\square$
\end{flushright}

\subsection*{2.9} In an extensive (or just Boolean)category, one has Sum $
\subseteq $ Mono, hence Mono$^{\square} \subseteq $ Sum$^{\square} =
$ SplitEpi. But in the presence of enough injectives,
Mono$^{\square}$ = SplitEpi only if Sum $ = $ Mono, a condition that
rarely holds even in a presheaf category: $\Set^{\C^{\text{op}}}$
satisfies Sum $ = \text{Mono}$ if, and only if, $\C$ is an
equivalence relation. For $\C = \{\xymatrix{\cdot \ar@<2pt>[r] \ar
@<-2pt>[r] & \cdot}\}$, so that $\Set^{\C^{\text{op}}}$is the
category of (directed multi-)graphs, with the Axiom of Choice
granted, Mono$^{\square}$ contais precisely the full morphisms that
are surjective on vertices; here a morphism $f: G \to H$ of graphs
is \emph{full} if every edge $f(a) \to f(b)$ in $H$ is the $f$-image
of an edge $a \to b$ in $G$.

\subsection*{2.10} For a wfs $(\E, \M)$ in a
category $\C$ with terminal object $1$, the full subcategory
\begin{equation*}
    \F(\M) := \{B \in \text{ob}\C \;|\; (B \to 1) \in \M\}
\end{equation*}
is weakly reflective in $\C$, in fact weakly $\E-$reflective, with
a weak reflection $\rho_A \in \E$ of an object $A$ being obtained
by an $(\E, \M)$-factorization of $A \to 1$:
\begin{equation*}
  \xymatrix{
    A \ar[r]^{\rho_A} & RA \ar[r]^{\M} & 1.\\
  }
\end{equation*}
If $(\E, \M)$ is an fs, $\F(\M)$ is $\E$--reflective in $\C$.
\subsection*{2.11 Remark} Weak factorization systems are abundant in homotopy theory.
In fact, a \emph{Quillen model category} $\C$ is defined as a
complete and cocomplete category together with three classes of
morphisms $\E$ (\emph{cofibrations}), $\M$ (\emph{fibrations}) and
$\W$ (\emph{weak equivalences}) such that $\W$ has the 3--for--2
property, is closed under retracts in $\C^2$ and $(\E,\M_0)$,
$(\E_0,\M)$ are weak factorization systems where
\[
\M_0=\M\cap\W, \quad\E_0=\E\cap\W
\]
denote the classes of trivial fibrations and cofibrations,
respectively. The \emph{3--for--2 property} means that whenever two
of the morphisms $gf$, $f$ and $g$ lie in $\W$, the third one lies
also in $\W$.

Objects of the weakly reflective subcategory $\F(\M)$ are called
\emph{fibrant}. Dually, when $\C$ has an initial object 0, there is
a weakly coreflective subcategory
\[
\T(\E) = \{A \in \ob \C\;|\; (0 \to A) \in \E\}
\]
of \emph{cofibrant} objects.

\section{Reflective factorization systems and prefibrations}

\subsection*{3.1} For a factorization system $(\E, \M)$ in a
category $\C$ with terminal object $1$, the $\E$--reflective full
subcategory $\F(\M)$ of 2.10 is even \emph{firmly $\E$-reflective},
in the sense that any morphism $A \to B$ in $\E$ with $\B \in
\F(\M)$ serves as a reflection of the object $A$ into $\F(\M)$. Such
reflective subcategories are easily characterized:

\subsection*{3.2 Proposition} \emph{ For a factorization system
$(\E, \M)$ and an $\E$-reflective subcategory $\F$ of $\C$, the
following conditions are equivalent:}
\begin{enumerate}
    \item[(i)] $\F = \F(\M)$,
    \item[(ii)]\emph{$\F$ is firmly $\E$-reflective in $\C$,}
    \item[(iii)]$\E \subseteq R^{-1}(\text{Iso}\C)$.
\end{enumerate}
\emph{If these conditions hold, one has $\E = R^{-1}(\iso\C)$ if,
and only
    if, $\E$ satisfies (in addition to} (3) \emph{of} 2.2 \emph{) the cancellation property}

    (4)\emph{ if $gf \in \E$ and $g \in \E$, then $f \in \E$.}

\begin{proof}
(i)$\implies$(ii): see 3.1. (ii) $\implies$ (iii): Considering the
$\rho$-naturality diagram for $e:A \to B$ in $\E$,
\begin{equation*}
  \xymatrix{
    A \ar[r]^{e} \ar[d]_{\rho_A} & B \ar[d]^{\rho_B}\\
    RA \ar[r]^{Re} & RB
  }
\end{equation*}
we see that $\rho_Be$ serves as a reflection for $A$, by
hypothesis, so that $Re$ must be an isomorphism. (iii) $\implies$
(i): For $B \in \F$, consider the $(\E, \M)$-factorization
\begin{equation*}
  \xymatrix{
    B \ar[r]^{e} & C \ar[r]^{m} & 1.\\
  }
\end{equation*}
Since $1 \in \F$ and $m \in \M$, also $\C$ lies in the
$\E$-reflective subcategory $\F$. Hence $e \cong Re$ is an
isomorphism, by hypothesis, and $(B \to 1)\in \M$ follows.
Conversely, having $(B \to 1) \in \M$, $1 \in \F$ implies $B \in
\F$, as above. $\E = R^{-1}(\text{Iso}\C)$ trivially implies (4),
and (4) implies $R^{-1}(\text{Iso}\C) \subseteq \E$, by inspection
of the $\rho$-naturality diagram above.
\end{proof}

We adopt the terminology of [CHK] and call an fs $(\E, \M)$ in any
category $\C$ \emph{reflective} if (4) holds. Since $\E$ is always
closed under composition and satisfies (3) of Section 2, we see that
an fs $(\E, \M)$ is a reflective fs if, and only if, $\E$ satisfies
the 3--for--2 property, granted the existence of cokernelpairs in
$\C$ (see 1.3).

A reflective fs $(\E_0, \M)$ makes $\C$ a Quillen model category,
with $\W = \E_0,\; \E = \C$ and $\M_0 = \iso\C$. The corresponding
homotopy category $\C[\W^{-1}]$ is $\F$.

\subsection*{3.3} A reflective fs $(\E, \M)$ in a category with
terminal object depends only on the reflective subcategory $\F(\M)$,
since $\E = R^{-1}(\text{Iso}\C)$ and $\M = \E^{\bot}$. Conversely,
given any reflective subcategory $\F$ of $\C$ with reflector $R$ and
reflection morphism $\rho:1 \to R$, one may ask when is $\E :=
R^{-1}(\text{Iso}\C)$ part of (a necessarily reflective) fs. This
question is discussed in general in [CHK], [JT1]. Here we are
primarily interested in the case when, moreover, $\E$'s
factorization partner $\M = \E^{\bot}$ can be presented as
\begin{equation*}
    \M = \text{Cart}(R, \rho),
\end{equation*}
where Cart$(R,\rho)$ is the class of $\rho$-\emph{cartesian
morphisms}, i.e., of those morphisms whose $\rho$-naturality
diagram is a pullback.

\subsection*{3.4 Proposition} \emph{For a reflective subcategory
$\F$ of the finitely complete category $\C$ with reflection}
$\rho: 1 \to R$, $(\E, \M) = \big(R^{-1}(\text{Iso}\C),
\text{Cart}(R,\rho)\big)$ \emph{is a factorization system of $\C$
if, and only if, for every morphism $f: A \to B$, the induced
morphism $e = (f, \rho_A):A \to B \times_{RB}RA$ into the pullback
of $Rf$ along $\rho_B$ lies in $\E$. In this case, $\F = \F(\M)$.}

\begin{proof} See Theorem 4.1 of [CHK] or Theorem 2.7 of [JT1].
\end{proof}

Adopting again the terminology used in [CHK], we call a reflective
factorization system $(\E,\M)$ \emph{simple} if $\M = \text{Cart}(R,
\rho)$, that is: if the reflective subcategory $\F = \F(\M)$
satisfies the equivalent conditions of Proposition 3.4. We also make
use of Theorem 4.3 of [CHK]:

\subsection*{3.5 Proposition} \emph{For a reflective fs $(\E, \M)$
of a finitely complete category $\C$, in the notation of} 3.3
\emph{the following conditions are equivalent, and they imply
simplicity of $(\E,\M)$:
\begin{enumerate}
    \item[(i)] $\E$ is stable under pullback along morphisms in
    $\M$;
    \item[(ii)]$R$ preserves pullbacks of morphisms in $\M$ along
    any other morphisms;
    \item[(iii)] the pullback of a reflection $\rho_A: A \to RA$
    along a morphism in $\F$ is a reflection morphism.
\end{enumerate}
} Reflective factorization systems $(\E, \M)$ satisfying these
equivalent conditions are called \emph{semi--left exact}. The
reflective subcategory $\F$ is a \emph{semilocalization} of $\C$ if
property (iii) holds; equivalently, if the associated reflective fs
is semi--left exact. A reflective fs need not be simple, and a
simple fs need not be semi--left exact(see [CHK]).

\subsection*{3.6} A Quillen model category $\C$ is called {\it right proper} if every pullback
of a weak equivalence along a fibration is a weak equivalence (see
[H]). Since each weak equivalence $w$ has a factorization $w=w_2w_1$
where $w_1$ is a trivial cofibration and $w_2$ a trivial fibration,
and since (trivial) fibrations are stable under pullback, $\C$
\emph{is right proper if, and only if, trivial cofibrations are
stable under pullback along fibrations, that is: if the  wfs
$(\E_0,\M)$ of} 2.10 \emph{has property} 3.5(i). Hence, a semi--left
exact reflective fs $(\E_0, \M)$ makes $\C$ a right proper Quillen
model category with $(\E, \M_0) = (\C, \iso\C)$ and $\W = \E_0$.

\subsection*{3.7}
Simple and semi-left exact reflective factorization systems occur
most naturally in the context of fibrations. Hence, recall that a
functor $P:\C \to \B$ is a \emph{(quasi-)fibration} if the induced
functors
\[
    P_C: \C/C \to \B/PC
\]
have full and faithful right adjoints, for all $\C \in \ob\C$. Let
us call $P$ a \emph{prefibration} if, for all $C$, there is an
adjunction
\[
\xymatrix@C=15pt{P_C\ar@{-|}[r]^-{\eta}_-{\varepsilon}&I_C}\ .
\]
whose induced monad is idempotent. (Janelidze's notion of admissible
reflective subcategory $\B$ of $\C$ asks the right adjoints $I_C$ to
be full and faithful, so that each $P_C$ is a fibration, in
particular a prefibration; see [J], [CJKP].) With the notation
\[
    I_C: (g:B \to PC) \mapsto (v_g:g^{*}C \to C)
\]
we can state right adjointness of $P_C$ more explicitly, as follows:
for every morphism $g: B \to PC$ in $\B$ one has a commutative
diagram
\[
    \xymatrix{
        P(g^{*}C) \ar[r]^{Pv_g} \ar[d]_{\varepsilon_g} & PC
        \ar[d]^1 \\
        B \ar[r]^g & PC\\
    }
\]
in $\B$, and whenever
\[
    \xymatrix{
        PA \ar[r]^{Pf} \ar[d]_u & PC \ar[d]^1 \\
        B \ar[r]^g & PC\\
    }
\]
commutes in $\B$ (with $f:A \to C$ in $\C$), then there is a
unique morphism $t:A \to g^{*}C$ in $\C$ with $v_gt =f$ and
$\varepsilon_g\cdot Pt = u$. \\
If $u =1$, then $t = \eta_f$, and we obtain the factorization
\[
\xymatrix {
    & (Pf)^{*}C \ar[dr]^{v_{Pf}} & \\
    A \ar[ur]^{\eta_f} \ar[rr]^f & & C\\
    }
\]
and the idempotency condition amounts to the requirement that
$P\eta_f = \varepsilon_{Pf}^{-1}$ be an isomorphism. One then has
$v_{Pf} \in \cart P$, with
\[
    \cart P = \{f \;|\; \eta_f \text{ iso}\}.
\]
(As we will see shortly, there is no clash with the notation used in
3.3.) In fact, $\big(P^{-1}(\iso \B), \cart P\big)$ is a
factorization system of $\C$, and it is trivially reflective.

Let us now \emph{assume that $P$ preserves the terminal object $1$
of $\C$}. Then
\[
    \F(\cart P) = \{A \;|\; A \to 1 \; P\text{--cartesian}\}
\]
contains precisely the $P$--\emph{indiscrete} objects of $\C$,
e.g. those $A \in \ob\C$ for which every $h:PD \to PA$ in $\B$
(with $D \in \ob \C$) can be written uniquely as $h = Pd$, with
$d: D \to A$ in $\C$. If we denote the adjunction
\[
   \xymatrix{ \C \simeq  \C/1 \ar@<1ex>[r]^{P_1} & \B/P1 \simeq \B \ar@<1ex>[l]^{I_1}}
\]
simply by $\xymatrix@C=15pt{P\ar@{-|}[r]^-{\eta}_-{\varepsilon}&I:
\B \to \C}$, then
\[
    \F(\cart P) = \{ A \;|\;\eta_A \; \text{ iso}\}
\]
is the reflective subcategory of $\C$ fixed by the adjunction
$\xymatrix@C=15pt{P\ar@{-|}[r]&I}$. Hence its reflector $R$ (as an
endofunctor of $\C$) is $IP$, with reflection morphism $\eta$.

A routine exercise shows
\[
    P^{-1}(\iso B) = R^{-1} ( \iso \C), \cart P = \cart (R, \eta).
\]
In particular, \emph{ the fs $\big(P^{-1}(\iso \B), \cart P\big)$
given by a prefibration $P$ with $P1 \cong 1$ is simple}. An easy
calculation shows also that $P^{-1}(\iso \B)$ is stable under
pullback along morphisms in $\cart P$ when $P$ preserves such
pullbacks. Consequently, \emph{for $\C$ finitely complete and with
the prefibration $P$ preserving pullbacks of $\cart P$--morphisms
and the terminal object, the fs is actually semi-left exact}.

\subsection*{3.8} Conversely to 3.7, let us show that any simple
reflective fs $(\E, \M)$ of a finitely complete category $\C$ is
induced by a prefibration $P$ with $P1 \cong 1$. More precisely, we
show that the restriction $\C \to \F(\M)$ of the reflector $R$
(notation as in 3.3) is a prefibration. To this end, for $g: B \to
RC$ with $B \in \F(\M)$ we form the (outer) pullback diagram
\[
    \xymatrix{
        B \times_{RC} C \ar[rr]^{v_g} \ar[dd]_p \ar[dr]^{\rho_{B\times_{RC}C}} & & C
        \ar[d]^{\rho_C}\\
        & R(B \times_{RC}C) \ar[r]_{Rv_g} \ar[dl]^{\varepsilon_g}
        & RC \ar[d]^1\\
        B \ar[rr]_g & & RC
    }
\]
The pullback projection $p$ factors through $R(B \times_{RC}C)$ by a
unique morphism $\varepsilon_g$ since $B \in \F(\M)$. To verify the
required universal property, consider $f:A \to C$ and $u:RA \to B$
with $Rf = gu$. Since
\[
    gu\rho_A = Rf\cdot\rho_A = \rho_Cf,
\]
there is a unique morphism $t:A \to B\times_{RC}C$ with $pt
=u\rho_A$, $v_gt = f$. From
\[
    \varepsilon_g\cdot Rt\cdot\rho_A =
    \varepsilon\rho_{B\times_{RC}C}t = pt = u\rho_A
\]
one obtains $\varepsilon_g\cdot Rt = u$, as required. Since,
conversely, $\varepsilon_g\cdot Rt = u$ implies $pt = u\rho_A$, we
have shown right adjointness of $R_C$. Furthermore, when $u =1$, the
pullback diagram above can simply be taken to be the
$\rho$--naturality diagram of $f$, by simplicity of $(\E, \M)$.
Hence, $A \cong B \times_{RC}C$ and $p \cong \rho_A$, so that
$\varepsilon_{Pf}$ is an isomorphism, and this shows the required
idempotency. Consequently, the reflector of $\F(\M)$ is a
prefibration, and since $\E = R^{-1}(\iso \C)$, the induced
factorization system must be the given fs $(\E, \M)$. By 3.5, the
system is semi-left exact precisely when the reflector preserves
pullbacks of morphisms in $\M$. Hence, with 3.7 we proved here:

\subsection*{3.9 Theorem}\emph{In a finitely complete category $\C$,
$(\E, \M)$ is a simple reflective factorization system of $\C$ if,
and only if, there exists a prefibration $P:\C \to \B$ preserving
the terminal object with
\[
    \E = P^{-1}(\iso B),\; \M = \cart P.
\]
$(\E,\M)$ is semi--left exact precisely when $P$ can be chosen to
preserve every pullback along a $P$--cartesian morphism}.
\begin{flushright}
$\Box$
\end{flushright}

\section{Torsion Theories}

\subsection*{4.1} Let $(\E, \M)$ be a reflective fs in a category
$\C$ with zero object $0 = 1$. (There is no further assumption on
$\C$ until 4.6.) Then we have not only the $\E$-reflective
subcategory $\F = \F(\M)$ with reflection $\rho: 1 \to R$, but also
the $\M$-coreflective subcategory $\T = \T(\E)$ (see 2.11), whose
coreflections $\sigma_B:S_B \to B$ are obtained by $(\E,
\M)$-factoring $0 \to B$, for all $B$ in $\C$. Let us first clarify
how $\T$ and $\F$ are related.

\subsection*{4.2 Proposition} \emph{In the setting of} 4.1, \emph{the
following assertions are equivalent for an object $A$ in $\C$:
\begin{enumerate}
    \item[(i)] $A \in \T$;
    \item[(ii)] $\C(A,B) = \{0\}$, for all $B \in \F$;
    \item[(iii)] $RA  \cong 0$.
\end{enumerate}
} \begin{proof} (i)$\implies$(ii) follows from $(0 \to A) \bot (B
\to 0)$. (ii)$\implies$(iii): Since $RA \in \F$, one has $\rho_A =
0$ and obtains $1_{RA} = 0$ from $\rho_A \bot (RA \to 0)$.
(iii)$\implies$(i): Since $RA \cong 0$, one has $(A \to 0) \in \E$,
and this implies $(0 \to A) \in \E$ by (4) of 3.2, hence $A \in \T$.
\end{proof}

Dualizing Propositions 3.2 and 4.2 we obtain:

\subsection*{4.3 Corollary} \emph {In the setting of }4.1, $\M =
S^{-1}(\text{Iso}\C)$\emph{ if, and only if, $\M$ satisfies the
cancellation property}:
\begin{enumerate}
    \item[(4$^{\text{op}}$)] \emph {if $gf \in \M$ and $f \in \M$,
    then $g \in \M$.}
\end{enumerate}
\emph{In this case,}
\begin{equation*}
    \F = \{B \in \text{ob}\C\;|\; SB \cong 0\}=\{B \;|\; \C(A,B)
    = \{0\} \text{ for all } A \in \T\}.
\end{equation*}
Factorization systems $(\E, \M)$ satisfying (4$^{\text{op}}$) are
called \emph{coreflective}.

\subsection*{4.4 Definitions and Summary} A \emph{torsion theory} in a
category $\C$ is a reflective and coreflective factorization system
$(\E, \M)$ of $\C$, i.e., a fs of $\C$ in which both classes satisfy
the 3--for--2 property. If $\C$ has kernelpairs or cokernelpairs, it
actually suffices to assume that $(\E, \M)$ be a wfs in this
definition (see 2.7, 2.8). If $\C$ has a zero object, then $\T =
\T(\E)$ is the \emph{torsion subcategory} and $\F = \F(\M)$ the
\emph{torsion-free subcategory} associated with the theory. For an
object $C$, the coreflection $\sigma_C$ into $\T$ and the reflection
$\rho_C$ into $\F$ are obtained by $(\E, \M)$-factoring $0 \to C$
and $C \to 0$, respectively as in
\begin{equation*}
  \xymatrix{
   0 \ar[r] & SC \ar[r]^{\sigma_C} & C \ar[r]^{\rho_C} & RC \ar[r] & 0.\\
  }
\end{equation*}
$R$ and $S$ determine all $\E, \M, \T, \F$, via
\begin{equation*}
    \begin{array}{cc}
      \E = R^{-1}(\text{Iso}\C) =\; ^{\bot}\M, & \M = S^{-1}(\text{Iso}\C) = \E^{\bot}, \\
      \T = R^{-1}(\{0\}) = \F^{\leftarrow}, & \F = S^{-1}(\{0\}) = \T^{\rightarrow}, \\
    \end{array}
\end{equation*}
with $\F^{\leftarrow} := \big\{A\;|\; \forall B \in \F \big(
\C(A,B) = \{0\}\big)\big\}, \T^{\rightarrow} := \big\{B\;|\;
\forall A\in \T \big(\C(A,B) = \{0\}\big)\big\}$. Furthermore, if
$\C$ has pullbacks and $\E$ is stable under pullback along
morphisms in $\M$, i.e., if the torsion theory is
\emph{semi-left-exact} and, hence, \emph{simple}, then an $(\E,
\M)$-factorization of $f:A \to B$ can be presented as
\begin{equation*}
  \xymatrix{
   & RA\times_{RB} B  \ar[dr]^{\pi_2} & \\
   A \ar[rr]_f \ar[ur]^{(\rho_A,f)} & & B
  }
\end{equation*}
where $\pi_2$ is the pullback of $Rf$ along $\rho_B$. In this
case, $\M = \text{Cart}(R, \rho)$. We note that without the
hypothesis of semi-left-exactness or simplicity, one still has:
\begin{equation*}
    f \in \E \iff \pi_2 \text{ iso,  } f \in \M \iff (\rho_A, f) \in \M.
\end{equation*}
The condition dual to semi-left-exactness is called
\emph{semi-right-exactness}, and it yields $\E =
\text{Cocart}(S,\sigma)$, along with an alternative presentation
of the $(\E, \M)$-factorization of $f$:
\begin{equation*}
  \xymatrix{
   & A+_{SA} SB  \ar[dr]^{(f, \sigma_B)} & \\
   A \ar[ur]^{\kappa_1} \ar[rr]_f & & B
  }
\end{equation*}
where $\kappa_1$ is the pushout of $Sf$ along $\sigma_A$.

In a category $\C$ with zero object, let 0Ker be the class of
morphisms whose kernel is $0$, and 0Coker the class of morphisms
with zero cokernel. Note that $ \text{Mono} \subseteq \text{0Ker}$
and $\text{Epi} \subseteq \text{0Coker}$.

\subsection*{4.5 Proposition} \emph{In a category $\C$ with
$0$, any pair of full subcategories $\T = \F^{\leftarrow}$ and $\F =
\T^{\rightarrow}$ satisfies the following properties, for any
morphisms $k:A \to B,  p:B \to C$ in $\C$}.
\begin{enumerate}
    \item[(1)] \emph{for} $k \in \text{0Ker}, B \in \F$
    \emph{implies} $A \in \F$;
    \item[(2)] \emph{for} $p \in \text{0Coker}, B \in \T$
    \emph{implies} $C \in \T$;
    \item[(3)] \emph{for $k$ the kernel of $p$, $A,C \in \F$
    imply $B \in \F$};
    \item[(4)] \emph{for $p$ the cokernel of $k$, $A,C \in \T$
    imply $B \in \T$.}
\end{enumerate}

\begin{proof}
(3) implies (1), and (2), (4) are dual to (1), (3), respectively.
Hence, if suffices to prove (3): any morphism $f:T \to B$ with $T
\in \T$ satisfies $pf = 0$. Hence, it factors through $k$, by a
morphism $T \to A$, which must be 0, so that also $f = 0$.
\end{proof}

\subsection*{4.6} We call a full subcategory $\F$ \emph{closed
under left-extensions} in $\C$ if it satisfies (3) of 4.5. If $\C$
has (NormEpi, 0Ker)-factorizations, with NormEpi the class of normal
epimorphisms (i.e. of morphisms that appear as cokernels), and if
$\F$ satisfies property (1) of 4.5, then the morphism $p$ in (3) may
be taken to be the cokernel of $k$, so that closure under
left-extensions amounts to the selfdual property of being
\emph{closed under extensions}. Note that $\C$ has (NormEpi,
0Ker)-factorization if $\C$ has kernels and cokernels (of kernels),
and if pullbacks of normal epimorphisms along normal monomorphisms
have cokernel 0 (see Prop. 2.1 of [CDT]). From 4.5 (1), (2) one
obtains:

\subsection*{4.7 Corollary} \emph{The reflection morphisms of the
torsion-free subcategory of a torsion theory in a pointed category
with} (NormEpi, 0Ker)-\emph{factorization are normal epimorphisms.
Dually, if there are} (0Coker, NormMono)--\emph{facto-rizations,
then the coreflection morphisms of the torsion subcategory are
normal monomorphisms.}

\subsection*{4.8} In a pointed category with kernels and
cokernels, let $(\E, \M)$ be a torsion theory. With the notation of
4.4, let $\kappa_C = \text{ker}\rho_C$ and $\pi_C =
\text{coker}\sigma_C$. If, as in 4.7, $\rho_C$ is a normal
epimorphism and $\sigma_C$ a normal monomorphism, so that $\rho_C =
\text{coker}\kappa_C$ and $\sigma_C = \text{ker}\pi_C$, we obtain
induced morphisms $\alpha_C$ and $\beta_C$ that, in the next
diagram, make squares 1, 2, 3 pullbacks and squares 2, 3, 4
pushouts:
\begin{equation*}
  \xymatrix{
    SC \ar[r]^{1} \ar[d]_{\alpha_C} \ar@{}[dr]|{\framebox{1}} & SC \ar[r] \ar[d]^{\sigma_C} \ar @{} [dr] |{\framebox{2}} & 0 \ar[d] \\
    KC \ar[r]^{\kappa_C} \ar[d] \ar@{}[dr]|{\framebox{3}} & C \ar[r]^{\pi_C} \ar[d]^{\rho_C} \ar @{} [dr] |{\framebox{4}} & QC \ar[d]^{\beta_C} \\
    0 \ar[r] &  RC \ar[r]_{1} & RC
  }
\end{equation*}
Since $\rho_C \in \E$, also $\beta_C \in \E$ (since $\E$ is pushout
stable), whence $\pi_C \in \E$ (by the 3-for-2 property) and
$R\pi_C$ iso. But since $\rho_{RC}$ is iso, this means that
$\beta_C$ may be replaced by $\rho_{QC}$. Likewise, replacing
$\alpha_C$ by $\sigma_{QC}$, we can redraw the above diagram as:
\begin{equation*}
  \xymatrix{
    SKC \ar[r]^{S\kappa_C}_{\sim} \ar[d]_{\sigma_{KC}} \ar@{}[dr]|{\framebox{1}} & SC \ar[r] \ar[d]^{\sigma_C} \ar @{} [dr] |{\framebox{2}} & 0 \ar[d] \\
    KC \ar[r]^{\kappa_C} \ar[d] \ar@{}[dr]|{\framebox{3}} & C \ar[r]^{\pi_C} \ar[d]^{\rho_C} \ar @{} [dr] |{\framebox{4}} & QC \ar[d]^{\rho_{QC}} \\
    0 \ar[r] &  RC \ar[r]_{R\pi_C} & RQC
  }
\end{equation*}
The endofunctors $K$ and $Q$ behave just like $S$ and $R$ when we
want to describe the subcategories $\T$ and $\F$:

\subsection*{4.9 Proposition} \emph{Under the hypothesis of} 4.8,
\emph{ for every object $C$ one has the following equivalences}:
\begin{equation*}
    \begin{array}{c}
      C \in \F(\M) \iff KC \in \F(\M) \iff KC = 0, \\
      C \in \T(\E) \iff QC \in \T(\E) \iff QC = 0. \\
    \end{array}
\end{equation*}

\begin{proof} Since $\kappa_C = \text{ker}\rho_C$ and $RC \in \F$, one has
$(C \in \F \iff KC \in \F)$ by Prop. 4.5. Furthermore, $(C \in \F
\iff \rho_C$ iso $\iff \kappa_C = 0 \iff KC = 0)$. The rest follows
dually. \end{proof}

The normal monomorphism $\alpha_C \cong \sigma_{KC}$ and the
normal epimorphism $\beta_C \cong \rho_{QC}$ measure the
``distance'' from $\kappa_C$ to the coreflection $\sigma_C$ and
from $\pi_C$ to the reflection $\pi_C$, respectively. The
following Theorem indicates when that ``distance'' is zero:

\subsection*{4.10 Theorem} \emph{Under the hypothesis of} 4.8, \emph{the
following conditions are equivalent for every object $C$}:
\begin{enumerate}
    \item[(i)] $\pi_C\cdot\kappa_C = 0$;
    \item[(ii)] $\ker\rho_{QC} = 0$;
    \item[(iii)] $\pi_{QC}$ \emph{is an isomorphism};
    \item[(iv)] $QC \in \F(\M)$;
    \item[(v)] $\text{coker}\sigma_{KC} = 0$;
    \item[(vi)] $\kappa_{KC}$ \emph{is an isomorphism};
    \item[(vii)] $KC \in \T(\E)$;
    \item[(viii)]$(0 \to QC) \in \M$;
    \item[(ix)] $(KC \to 0) \in \E$.
\end{enumerate}
\emph{All conditions are satisfied when $(\E, \M)$ is simple} (see
3.4).

\begin{proof} Since $\rho_{QC} \cdot \pi_C \cdot \kappa_C = 0$, (i) $\iff$
(ii) is obvious. (iv) implies $\rho_{QC}$ iso, hence (ii), and
also (iii), since
\begin{equation*}
    \rho_{QQC} \cdot \pi_{QC} = R\pi_{QC} \cdot \rho_{QC},
\end{equation*}
with $R\pi_{QC}$ iso. Conversely, (ii) $\implies$ (iv) holds since
$\rho_{QC}$ is a normal epimorphism, and (iii) $\implies$ (iv)
holds since $\pi_{QC} = \text{coker}\sigma_{QC}$ iso means $SQC =
0$, hence $QC \in \F(\M)$. Consequently, we have (i) $\iff$ (ii)
$\iff$ (iii) $\iff$ (iv), and (i) $\iff$ (v) $\iff$ (vi) $\iff$
(vii) follows dually. Since
\begin{equation*}
  \xymatrix{
   KC \ar[r]^{\rho_{KC}} & RKC \ar[r] & 0 \\
  }
\end{equation*}
is the $(\E, \M)$-factorization system of $KC \to 0$, one has $RKC
\to 0$ iso, if, and only if, $KC \to 0$ lies in $\E$. This shows
(vii) $\iff$ (ix), and (iv) $\iff$ (viii) follows dually. Finally,
assume $(\E, \M)$ to be simple and consider the commutative
diagram
\begin{equation*}
  \xymatrix{
    KC \ar[r]^{1} \ar[d]_{\rho_{KC}} & KC \ar[r]^{\kappa_C} \ar[d] \ar @{} [dr] |{\framebox{3}} & C \ar[d]^{\rho_{C}} \\
    RKC \ar[r] &  0 \ar[r] & RC
  }
\end{equation*}
Since $0 = \rho_C \kappa_C = R\kappa_C \cdot \rho_{KC}$ with
$\rho_{KC}$ epi, the buttom row is $R\kappa_C$, and diagram 1 of 4.9
shows that $\kappa_C$ lies in $\M$, since $(\E,\M)$ is coreflective.
Hence, the whole diagram is a pullback, by simplicity of $(\E,\M)$,
and therefore also its left square: $KC \cong KC \times RKC$. Now
the morphism $t = \langle 0:RKC \to KC, 1_{RKC}\rangle$ shows that
$\rho_{KC}$ must be 0, which means $RKC = 0$ and, hence, $KC \in
\T$.
\end{proof}

\subsection*{4.11 Remarks}
(1) Following the terminology of [CHK] we call a torsion theory
\emph{normal} if the equivalent conditions of 4.10 hold. Hence
\emph{every simple torsion theory is normal}, provided that $\C$
satisfies the hypothesis of 3.8. Moreover, square 3 of 4.8 and
condition (ix) of 4.10 show that $(\E, \M)$ is normal if, and only
if, $\E$ satisfies a very particular pullback-stability condition.
No failure of this condition is known since the following open
problem of [CHK] remains unsolved: \emph{is there a non-normal
torsion theory?}

(2) The advantage of our definition of torsion theory is that we
do not need to assume the existence of kernels and cokernels in
$\C$. It applies, for example, to a triangulated category $\C$.
Such a category has only weak kernels and weak cokernels and our
definition precisely corresponds to torsion theories considered
there as pairs $\F$ and $\T$ of colocalizing and localizing
subcategories (see [HPS]).

It is also easy to express torsion theories in terms of
prefibrations, since Theorem 3.9 gives immediately:

\subsection*{4.12 Corollary} \emph{In a finitely complete category $\C$,
the class $\M$ belongs to a torsion theory $(\E, \M)$ if, and only
if, there is a prefibration $P:\C \to \B$ with $P1 \cong 1$ such
that $\M = \cart P$ has the 3--for--2 property. Dually, in a
finitely cocomplete category $\C$, the class $\E$ belongs to a
torsion theory $(\E, \M)$ if, and only if, there is a precofibration
$Q:\C \to \A$ with $Q0 \cong 0$ such that $\E = \mathrm{Cocart}Q$
has the 3--for--2 property}.
\begin{flushright}
$\Box$
\end{flushright}

\section{Characterization of normal torsion theories}

\subsection*{5.1} In a finitely complete category $\C$ with a zero
object and cokernels (of normal monomorphisms), we wish to compare
the notion of normal torsion theory (as presented in 4.4, 4.11) with
concepts considered previously, specifically with the more classical
notion used in [BG] and [CDT]. Hence here let us refer to a pair
$(\T, \F)$ of full replete subcategories of $\C$ satisfying
\begin{enumerate}
    \item $\C(A,B) = \{0\}$ for all $A \in \T$ and $B \in \F$ ,
    \item for every object $C$ of $\C$ there exists $\xymatrix{A \ar[r]^k &  C \ar[r]^q &
    B}$ with $A \in \T,\; B \in \F, \;k = \text{ker}q,\; q = \text{coker}k$.
\end{enumerate}
as a \emph{standard torsion theory} of $\C$; its torsion-free part
is necessarily normal-epireflective in $\C$. The main result of
[JT2] states that, when normal epimorphisms are stable under
pullback in $\C$, \emph{a normal-epireflective subcategory $\F$ is
part of a standard torsion theory if, and only if, $\F$ satisfies
the following equivalent conditions:}

\begin{enumerate}
    \item[(i)]  $\F$ \emph{is a semilocalization of $\C$ (see}
    3.5\emph{);}
    \item[(ii)] \emph{the reflector $\C \to \F$ is a (quasi)fibration (see} 3.7 \emph{);}
    \item[(iii)] \emph{$\F$ is closed under extensions, and the pushout of the
    kernel $\xymatrix{A \ar[r]^k &C}$ of $\rho_C$ along $\rho_A$
    is a normal monomorphism, for every} $C \in \text{ob}\C$ \emph{(with
    $\rho_C$ the $\F$-reflection of $C$)}.
\end{enumerate}
Recall that $\C$ is \emph{homological} [BB] if it is regular
[Ba] and protomodular [Bo]; here the latter property amounts
to: if in the commutative diagram
\begin{equation*}
    \xymatrix{
        \cdot \ar[d] \ar[r] & \cdot \ar[d]^p \ar[r] & \cdot \ar[d] \\
        \cdot \ar[r] & \cdot \ar[r] & \cdot\\
    }
\end{equation*}
with regular epimorphism $p$ the left and the whole rectangles are
pullbacks, so is the right one. In such categories one has
(NormEpi, 0Ker) = (RegEpi, Mono).

We are now ready to prove:

\subsection*{5.2 Theorem} \emph{Every standard torsion theory of $\C$
determines a simple reflective factorization system $(\E, \M)$ of
$\C$ with $\F(\M)$ normal-epireflective and $\T(\E)$
normal-monocoreflective. When $\C$ is homological, such
factorization systems are normal torsion theories. When both $\C$
and} $\C^{\text{op}}$ \emph{are homological, then normal torsion
theories correspond bijectively to standard torsion theories.}

\begin{proof}
Since a standard torsion theory $(\T, \F)$ is given by the
semilocalization $\F$, its reflective factorization system $(\E,
\M)$ is simple (see 3.4, 3.5), and one has $\T = \T(\E)$ (see 4.2).
This proves the first statement. For the second, let $(\E, \M)$ be a
simple reflective factorization system such that the reflections of
$\F(\M)$ are normal epimorphisms and the coreflections of $\T(\M)$
are normal monomorphisms. Simplicity means $\M = \text{Cart}(R,
\rho)$ by 3.4, and since the reflections $\rho_C\;(C \in
\text{ob}\C)$ are regular epimorphisms, protomodularity of $\C$
gives immediately that $\M$ satisfies the 3--for--2 property. Hence
$(\E, \M)$ is a torsion theory, and its normality follows from 4.10,
which is applicable since the assumptions of 4.8 are fulfilled, by
hypothesis. When both $\C$ and $\C^{\text{op}}$ are homological,
because of 4.7 we can apply 4.10 and obtain the last statement.
\end{proof}

\subsection*{5.3 Remarks} (1) As the proof of 5.2 shows, for the
bijective correspondence between normal torsion theories and
standard torsion theories, it suffices to have $\C$ homological with
(0Coker, NormMono)-factorizations. The latter condition is, of
course, still quite restrictive: even standard semi-abelian
categories (like the categories of groups or of commutative rings)
do not satisfy it. However, the type of categories that are both
homological and co-homological is very well studied. As George
Janelidze observed, these are precisely the "Raikov semi-abelian"
[Ra], [K] or "almost-abelian" [Ru] categories. In fact, in a pointed
protomodular category, the canonical morphism $A+B \to A \times B$
is an extremal epimorphism, hence it is an isomorphism when the
category is also co-protomodular. Since protomodular categories are
Mal'cev, co-protomodularity makes such categries additive. Hence,
\emph{the following conditions are equivalent for a category $\C$:}

\begin{enumerate}
    \item[(i)]  $\C$ \emph{is regular, coregular and additive;}
    \item[(ii)] $\C$ \emph{is homological and co-homological;}
    \item[(iii)] $\C$ \emph{is Raikov semi-abelian (= almost-abelian)}.
\end{enumerate}

Clearly, these conditions imply that $\C$ is homological with
(0Coker, NormMono)-factorizations, but we don't know whether these
properties are equivalent to (i)-(iii).

(2) Consider the additive homological category $\C$ of abelian
groups satisfying the implication $(4x = 0 \implies 2x = 0)$. As
shown in [JT2], the subcategory $\F$ of groups satisfying $2x = 0$
is closed under extensions and normal epireflective, but is not part
of a standard torsion theory. Its reflective factorization system is
not simple (likewise when one considers it not in $\C$ but in the
abelian category of all abelian groups, see [CHK]), and it is not a
normal torsion theory of $\C$. In fact, for $C = \mathbb{Z}$, the
diagram of 4.8 is as follows:
\begin{equation*}
\xymatrix {
    0 \ar[r] \ar[d]^{\sigma} & 0 \ar[r] \ar[d]^{\sigma} & 0 \ar[d] \\
    \mathbb{Z} \cong 2\mathbb{Z} \ar@{^{(}->}[r]^{\kappa} \ar[d]& \mathbb{Z}  \ar[d]_{\rho} \ar[r]^{\pi = 1}&
    \mathbb{Z} \ar[d]^{\rho}\\
    0 \ar[r] & \mathbb{Z}_2 \ar[r]^{1} & \mathbb{Z}_2
}
\end{equation*}
But we do not know whether $(\E, \M)$ is a torsion theory.

\subsection*{5.4} A standard torsion theory is called {\it hereditary} if $\T$ is closed under normal
subobjects, and it is {\it cohereditary} if $\F$ is closed under
normal quotients. While hereditary standard torsion theories are
of principal importance, coheredity is a very restrictive
property, as we show in the next proposition, which is well--known
in the case of groups (see [N]).

\subsection*{5.5 Proposition} \emph{Let $\C$ be a pointed variety of universal algebras where free algebras
are closed under normal subobjects. Then each standard
cohereditary torsion theory $(\T,\F)$ in $\C$ is trivial, i.e.,
$\T=\C$ or $\F=\C$.}

\begin{proof}
Assume $\F\neq\C$. Since $\F$ is closed under normal quotients,
there is a free algebra $V$ not belonging to $\F$. Hence, the
$\T$--coreflection of $V$ satisfies
\[
0\neq KV\in\T,
\]
and $KV$ is free (as a normal subobject of a free algebra)and
belongs to $\T$. Since $\T$ is closed under coproducts and
quotients, $\T=\C$ follows.
\end{proof}

\end{document}